# CONTOU-CARRÈRE SYMBOL VIA ITERATED INTEGRALS

ZHENBIN LUO
DEPARTMENT OF MATHEMATICS
BRANDEIS UNIVERSITY

ABSTRACT. This paper gives a new definition of the Contou-Carrère symbol in terms of an exponential of a Chen iterated integral and proves the corresponding reciprocity law.

## 0. INTRODUCTION

Given a discrete valuation $\nu$ on a field $F$, with valuation ring $A_\nu$, unique maximal ideal $p_\nu$ and residue field $k_\nu$, J. Milnor[M] defined the tame symbol $d_\nu : F^\times \times F^\times \to k_\nu^\times$ by

$$d_\nu(x,y) = (-1)^{\nu(x)\nu(y)} \frac{x^{\nu(y)}}{y^{\nu(x)}} \pmod{p_\nu}$$

This symbol have since been interpreted and generalized in several ways, algebraically and geometrically. Arbarrello, de Concini and V.G. Kac[A-D-K] defined the tame symbol of an algebraic curve as the commutator of a certain central extension of group. Latter, F. Pablos Romo[P1] gave an interpretation of this central extension in terms of determinants associated with infinite-dimensional vector subspaces. On the geometric side, P. Deligne[D] proved the following theorem by constructing from $f$ and $g$, a holomorphic line bundle $(f,g]$ with holomorphic connection on the complement of divisors of $f$ and $g$.

**Theorem 0.1.** *[D][W] Let $f,g$ be two meromorphic functions on a compact Riemann Surface $X$. Then $\prod_p \{f,g\}_p = 1$, where $\{f,g\}_p$ is the tame symbol $(-1)^{\nu(f)\nu(g)}[f^{\nu(g)}/g^{\nu(f)}](p)$ and $\nu(f)$ denotes the valuation at $p$ of the function $f$.*

In 1994, C. Contou-Carrère [C] defined a natural transformation greatly generalizing the tame symbol. In the case of an artinian local base ring $A$ with maximal ideal $m$, the symbol is given by the following formula. Let $f,g$ be two functions in $A((X))^\times$, they can be written uniquely as

$$f = a_0 X^{\nu(f)} \prod_{i=-\infty}^{+\infty}(1 - a_i X^i)$$

$$g = b_0 X^{\nu(g)} \prod_{i=-\infty}^{+\infty}(1 - b_i X^i)$$

where $(f, g)$ is the divisor of $g$ and $g$, $\nu_f, \nu_g \in \mathbb{Z}$, $a_i, b_i \in A$ for $i > 0$, $a_0, b_0 \in A^\times$, $a_i, b_i \in m$ for $i < 0$, and $a_i, b_i$ are zero when $i \ll 0$. Then the symbol is given by

$$< f, g >_{A((X))^\times} := (-1)^{\nu_f \nu_g} \frac{a_0^{\nu_g} \prod_{j=1}^\infty \prod_{k=1}^\infty (1 - a_j^{k/(j,k)} b_{-k}^{j/(j,k)})^{(j,k)}}{b_0^{\nu_f} \prod_{j=1}^\infty \prod_{k=1}^\infty (1 - a_{-j}^{k/(j,k)} b_k^{j/(j,k)})^{(j,k)}}$$

Since $a_{-i}, b_{-i}$ are zero when $i$ is large, the product is actually finite, hence the definition makes sense. This symbol is antisymmetric and, although not obvious, bimultiplicative and satisfies the Steinberg property. G. Anderson and F. Pablos Romo[A-P] gave an interpretation of the symbol as a commutator of liftings of $f$ and $g$ to a certain central extension of a group containing $A((X))^\times$.

In this paper we construct the Contou-Carrère symbol in terms of Chen iterated integrals and prove the corresponding reciprocity using a geometric property of iterated integrals. Consider a local artinian $\mathbb{C}-$ algebra $A$, finitely generated by nilpotent elements $\{a_1, a_2, \cdots, a_n\}$ as a $\mathbb{C}$ module. Let $\omega$ be a meromorphic 1-form on a compact Riemann surface $X$ with coefficients in $A$, in other words, it is of the form $\sum_{i=0}^n a_i \omega_i$, where $a_0 = 1$ and $\omega_i$'s are regular meromorphic 1-forms. Let $\gamma$ be a path on the Riemann surface. Define the integral of $f$ by extending the regular integral linearly:

$$\int_\gamma \omega = \sum_{i=0}^n a_i \int_\gamma \omega_i$$

In particular, we consider those 1-forms with logrithmetic poles, in other words, of the form $\frac{df}{f}$, where $f$ is a meromorphic function on $X$ with coefficients in $A$. Let $f$ and $g$ be two such functions, and $s$ be a zero or pole of $f$ or $g$ or both. Write $f$ and $g$ in Laurent series in powers of a uniformizor $x_s$, as elements in $A((x_s))^\times$. Thus, as mentioned above, we can express $f$ and $g$ locally as infinite products. Let $\sigma$ be a simple loop that starts and ends at some fixed point $P$, go around the divisor $s$ once in the counterclockwise direction but not any other divisors of $f$ or $g$. Then we have the following result:

**Theorem.**

$$\exp(\frac{1}{2\pi i} \int_\sigma \frac{df}{f} \circ \frac{dg}{g}) = (-1)^{\nu_f \nu_g} \frac{g(P)^{\nu_f} a_0^{\nu_g} \prod_{j=1}^\infty \prod_{k=1}^\infty (1 - a_j^{k/(j,k)} b_{-k}^{j/(j,k)})^{(j,k)}}{f(P)^{\nu_g} b_0^{\nu_f} \prod_{j=1}^\infty \prod_{k=1}^\infty (1 - a_{-j}^{k/(j,k)} b_k^{j/(j,k)})^{(j,k)}}$$

This formula differs from the Contou-Carrère symbol only by $g(P)^{\nu_f}/f(P)^{\nu_g}$, but since the point $P$ is fixed, we can rescale the functions $f$ and $g$ so that $f(P) = g(P) = 1$. Then we can get the usual Contou-Carrère symbol. This gives us a new interpretation of the Contou-Carrère in terms of exponential of an iterated integral, which is very convenient to use. Many properties of the symbol easily follows from the properties of iterated integrals. If we can consider the iterated integral over a loop that represents a relation of the fundamental group, by the homotopy invariance of iterated integral, the result is 0. Using this, we can reprove the reciprocity of Contou-Carrère easily. When the functions $f$ and $g$ are regular meromorphic functions, we get back the tame symbol and the Weil reciprocity.

# 1. Preliminaries of iterated integral

## 1.1. Definition of iterated integral.
For proofs of theorems of this section, see Chen[Ch] or Goncharov[G].

**Definition 1.1.** Let $\omega_1, \omega_2, \cdots, \omega_n$ be holomorphic 1-forms on a simply connected open subset $U$ of the complex plane $\mathbb{C}$. Let $\gamma : [0,1] \to U$ be a path. Then we call the integral

$$\int_\gamma \omega_1 \circ \cdots \circ \omega_n := \int \cdots \int_{0 \leq t_1 \leq \cdots \leq t_n \leq 1} \gamma^*\omega_1(t_1) \wedge \cdots \wedge \gamma^*\omega_n(t_n)$$

the iterated integral of the differential forms $\omega_1, \omega_2, \cdots, \omega_n$ over the path $\gamma$.

## 1.2. Homotopy invariance of iterated integrals. [H]

**Theorem 1.2.** Let $\omega_1, \cdots, \omega_n$ be holomorphic 1-forms on a simply connected open subset $U$ of the complex plane $\mathbb{C}$. Let $H : [0,1] \times [0,1] \to U$ be a homotopy, fixing the end points, of paths $\gamma_s : [0,1] \to U$ such that $\gamma_s(t) = H(s,t)$. Then

$$\int_{\gamma_s} \omega_1 \circ \cdots \circ \omega_n$$

is independent of $s$.

## 1.3. Shuffle relation.

**Theorem 1.3.** *[Shuffle relation]* Let $\omega_1, \cdots, \omega_n, \omega_{n+1}, \cdots, \omega_{n+n}$ be differential 1-forms, some of them could repeat. Let also $\gamma$ be a path that does not pass through any of the poles of the given differential forms. Denote by $Sh(m,n)$ the shuffles, which are permutations $\tau$ of the set $\{1, ..., m, m+1, ..., m+n\}$ such that $\tau(1) < \tau(2) < \cdots < \tau(m)$ and $\tau(m+1) < \tau(m+2) < \cdots < \tau(m+n)$. Then

$$\int_\gamma \omega_1 \circ \cdots \circ \omega_n \int_\gamma \omega_{n+1} \circ \cdots \circ \omega_{m+n} = \sum_{\tau \in Sh(m,n)} \int_\gamma \omega_{\tau(1)} \circ \omega_{\tau(2)} \cdots \circ \omega_{\tau(m+n)}.$$

## 1.4. Reversing the path.

**Lemma 1.4** (Reversing the path).

$$\int_\gamma \omega_1 \circ \omega_2 \circ \cdots \circ \omega_n = (-1)^n \int_{\gamma^{-1}} \omega_n \circ \omega_{n-1} \circ \cdots \circ \omega_1$$

## 1.5. Composition of paths.

**Theorem 1.5.** *[Composition of paths]* Let $\omega_1, \omega_2, \cdots, \omega_n$ be differential forms, some of them could repeat. Let $\gamma_1$ be a path that ends at $Q$ and $\gamma_2$ be a path that starts at $Q$. Then

$$\int_{\gamma_1\gamma_2} \omega_1 \circ \omega_2 \circ \cdots \circ \omega_n = \sum_{i=0}^n \int_{\gamma_1} \omega_1 \circ \omega_2 \circ \cdots \circ \omega_i \int_{\gamma_2} \omega_i \circ \omega_{i+1} \circ \cdots \circ \omega_n$$

## 2. Reciprocity law from a differential equation

Let $\omega_1, \omega_2, \cdots, \omega_n$ be meromorphic 1-forms on a compact Riemann surface $X$. Let $A_1, A_2, \cdots, A_n$ be non-commuting formal variables. Let $Y$ be the subset obtained from $X$ by removing the poles of $\omega_1 \circ \omega_2 \circ \cdots \circ \omega_n$. Consider the differential equation on $Y$:

$$dF = F(\sum_{i=1}^{n} A_i \omega_i).$$

Fix a point $p$ in $Y$. Let $\gamma : [0,1] \to Y$ be a piecewise smooth path starting at $p$ and ending at $z$. Think of the end point $z$ as a variable. Then the solution of the differential equation with initial condition $F(p) = 1$ is given by:

$$F_\gamma = 1 + \sum_{i=1}^{n} A_i \int_\gamma \omega_i + \sum_{1 \le i \le j \le n} A_i A_j \int_\gamma \omega_i \omega_j + \cdots$$

If $\gamma_1$ and $\gamma_2$ be two paths such that the end point of $\gamma_1$ is the beginning point of $\gamma_2$, then we have $F_{\gamma_1 \gamma_2} = F_{\gamma_1} F_{\gamma_2}$.

Consider loops in $Y$ that start at $P$ counterclockwise, bounding a disk in $X$ so that the disk contains only one poles of $\omega_1, \omega_2, \cdots, \omega_n$. We can choose these loops so that they do not intersect each other except at the point $P$. We can start with loops that bound a pole of $\omega_1$. Then continue with loops that bound a pole of $\omega_2$ and so on. Call these loops $\sigma_1, \sigma_2, \cdots, \sigma_N$. Let also $\alpha_i$ and $\beta_i$ be the loops in Y for $i = 1, ..., g$, where $g$ is the genus of $Y$ such that $\alpha_i$ and $\beta_i$ and $\sigma_1, \sigma_2, \cdots, \sigma_N$ generate $\pi_1(Y)$ and the only relation between them is

$$\sigma_1 \sigma_2 \cdots \sigma_N [\alpha_1, \beta_1] \cdots [\alpha_g, \beta_g] = 1.$$

If we multiply all the loops going in counterclockwise direction, then we will obtain a loop homotopic to the zero loop at $P$. This gives the global reciprocity law which can be formulate as follows if we view $F_{\sigma_i}$ as an element in the power series ring $\mathbb{C}[[A_1, A_2, \cdots, A_n]]$:

$$F_{\sigma_1} F_{\sigma_2} \cdots F_{\sigma_N} F_{[\alpha_1, \beta_1]} \cdots F_{[\alpha_g, \beta_g]} = 1$$

If we consider two differential forms $\omega_1, \omega_2$, then by comparing the coefficients of $A_1 A_2$ in the above equality, after some simplification we have the following equality:

$$\sum_{i=1}^{N} \int_{\sigma_i} \omega_1 \circ \omega_2 + \sum_{1 \le i \le j \le N} \int_{\sigma_i} \omega_1 \int_{\sigma_j} \omega_2 + \sum_{i=1}^{g} (\int_{\alpha_i} \omega_1 \int_{\beta_i} \omega_2 - \int_{\beta_i} \omega_1 \int_{\alpha_i} \omega_2) = 0$$

To see why this is true, first notice that $F_{[\alpha_i, \beta_i]}$ has no linear terms since the integral of a 1-form over such a commutator is zero. Hence the only contributors are quadratic terms from $F_{\sigma_i}$ and $F_{[\alpha_j, \beta_j]}$, and products of linear terms from $F_{\sigma_i}$ and $F_{\sigma_j}$. Using the product formula, it is easy to prove that the quadratic term of $F_{[\alpha_j, \beta_j]}$, ie. $\int_{[\alpha_j, \beta_j]} \omega_1 \circ \omega_2$ is equal to $\int_{\alpha_j} \omega_1 \int_{\beta_j} \omega_2 - \int_{\beta_j} \omega_1 \int_{\alpha_j} \omega_2$. For a more detailed proof, see [H].

If $\omega_1, \omega_2$ are differential forms of the third kind, for example, if $\omega_1 = \frac{df}{f}, \omega_2 = \frac{dg}{g}$, then the Weil reciprocity law can be easily reproved from this equality.

## 3. A NEW DEFINITION OF THE CONTOU-CARRÈRE SYMBOL

Let $X$ be a Riemann surface. By considering invertible rational functions with coefficient in a local artinian $\mathbb{C}$-algebra A, finitely generated as a $\mathbb{C}$-module, using the equality derived from section 2, we can obtain an equivalent definition of the Contou-Carrère symbol in terms of exponential of an integral. Based on this new definition, we can derive many properties of the symbol easily.

Let $A$ be a finite local $\mathbb{C}$-algebra with maximal ideal $m$. Let $f$ and $g$ be two invertible meromorphic functions with coefficient in $A$. In other words, $f$ and $g$ are functions of the form $\sum_{i=0}^{n} a_i f_i$, where all $f_i$ are regular meromorphic functions and $a_0 = 1, f_0 \neq 0, a_i \in m$ for $i$ not equal to 0. We say such a function is holomorphic if all $f_i$'s are holomorphic. Define the integral of such function by extending the usual integral linearly as mentioned in the introduction. It is easy to see that the sums of residues of such functions are still 0. Hence the reciprocity law for two differential forms derived in the previous section still holds for differential forms with coefficients in $A$. In particular, let us consider the case $\omega_1 = \dfrac{df}{f}, \omega_2 = \dfrac{dg}{g}$, forms with logarithmic poles. Let $s$ be a zero or pole of $f$ or $g$ or both. We can write $f$ and $g$ in Laurent series in powers of a uniformizor $x_s$, ie. $f$ and $g$ can be expressed locally as elements in $A((x_s))^\times$. Denote them by $f_s$ and $g_s$. By the product formula, we can write them as follows:

$$f_s = a_0 x^{\nu_f} \prod_{j=-N_f}^{-1} (1 - a_j x_s^j) \prod_{j=1}^{+\infty} (1 - a_j x_s^j)$$

$$g_s = b_0 x^{\nu_g} \prod_{j=-N_g}^{-1} (1 - b_j x_s^j) \prod_{j=1}^{+\infty} (1 - b_j x_s^j)$$

where $a_0$ and $b_0$ belong to $A^\times$, $a_j$ and $b_j$ belong to the maximal ideal $m$ for $j < 0$, and $a_j$ and $b_j$ belong to $A$ for $j > 0$. In the following, I will use $x$ instead of $x_s$ for short, and $f, g$ instead of $f_s, g_s$. Denote $\prod_{j=-N_f}^{-1}(1 - a_j x^j)$ by $f^-$, $\prod_{j=1}^{+\infty}(1 - a_j x^j)$ by $f^+$, and similarly for $g$. Let $\sigma$ be a loop that starts and ends at some point $P$, go around the divisor $s$ once in the counterclockwise direction but not around any other divisors of $f$ or $g$. Then we have the following result:

**Theorem.**

$$\exp(\frac{1}{2\pi i} \int_\sigma \frac{df}{f} \circ \frac{dg}{g}) = (-1)^{\nu_f \nu_g} \frac{g(P)^{\nu_f} a_0^{\nu_g} \prod_{j=1}^{\infty} \prod_{k=1}^{\infty} (1 - a_j^{k/(j,k)} b_{-k}^{j/(j,k)})^{(j,k)}}{f(P)^{\nu_g} b_0^{\nu_f} \prod_{j=1}^{\infty} \prod_{k=1}^{\infty} (1 - a_{-j}^{k/(j,k)} b_k^{j/(j,k)})^{(j,k)}}$$

**Remark 3.1.** *This gives us a new definition of the Contou-Carrère in terms of exponential of an iterated integral, which is very convenient to use.*

Denote the right hand side by $< f, g >_s$ for short. Let $\gamma_\varepsilon$ be a path that starts at $P$ and ends at some $Q_\varepsilon$ very close to 0, say of distance $\varepsilon$. Let $\sigma_\varepsilon$ be a loop that starts and ends at $Q_\varepsilon$ and goes around $s$ once in a counterclockwise direction along the circle of radius $\varepsilon$ and centered at $s$. Then we can deform $\sigma$ homotopically to $\gamma_\varepsilon \sigma_\varepsilon \gamma_\varepsilon^{-1}$. To prove the theorem, we need the following lemmas:

**Lemma 3.2.**
$$\int_{\sigma_\varepsilon} \frac{dx}{x} \circ \frac{dx}{x} \circ \cdots \circ \frac{dx}{x} = \frac{(2\pi i)^r}{r!}$$
where $\frac{dx}{x}$ is iterated $r$ times.

*Proof.* Take the parametrization of $\sigma_\varepsilon : x \to \varepsilon e^{2\pi i t}, 0 \le t \le 1$. Then by definition, we have
$$\int_{\sigma_\varepsilon} \frac{dx}{x} \circ \frac{dx}{x} \circ \cdots \circ \frac{dx}{x} = \int_0^1 \frac{dt}{t} \circ \cdots \circ \frac{dt}{t} = \frac{(2\pi i)^r}{r!}.$$

**Lemma 3.3.**
$$\int_{\sigma_\varepsilon} \frac{df}{f} = 2\pi i Res_{x=0} \frac{df_0}{f_0} = 2\pi i \nu(f)$$

*Proof.* Since $f(x) = \sum_{i=0}^n a_i f_i = f_0(1+F)$, where $F = f_0^{-1} \sum_{i=1}^n a_i f_i \in m(t)^\times$, we have
$$\frac{df}{f} = \frac{df_0}{f_0} + \frac{d(1+F)}{1+F}.$$

But
$$\frac{d(1+F)}{1+F} = d(\sum_{j=1}^N (-F)^j/j),$$
when integrated over a loop, this gives 0. Hence the lemma is proved.

If $a$ is an element of the maximal ideal $m$, then we define $\log(1-a)$ and $e^a$ as follows:
$$\log(1-a) := -\sum_{k=1}^\infty \frac{a^k}{k}, \quad e^a := \sum_0^\infty \frac{a^n}{n!}.$$

Since $m$ is nilpotent, the summation in the definitions are finite sums, so they are well-defined. Under these definitions, $\log(1-a)$ is an element of $m$ and $e^a$ is an element of $A^\times$. It is easy to see that we have the following relations as usual:
$$e^{\log(1-a)} = 1-a, \quad \log e^a = a.$$

**Lemma 3.4.**
$$\int_{\sigma_\varepsilon} \frac{dx}{x} \circ \frac{d(1-ax^n)}{1-ax^n} = 2\pi i \log(1-a\varepsilon^n)$$
where $n$ is an integer.

*Proof.* Take the parametrization of $\sigma_\varepsilon : x \to \varepsilon e^{2\pi i t}, 0 \leq t \leq 1$. Then by definition, we have

$$\begin{aligned}
\int_{\sigma_\varepsilon} \frac{dx}{x} \circ \frac{d(1-ax^n)}{1-ax^n} &= \int_0^1 \left(\int_0^{t_2} 2\pi i dt_1\right)\left(\sum_{k=0}^\infty (ax^n)^k(-2\pi i a n x^n)\right) dt_2 \\
&= -(2\pi i)^2 \sum_{k=0}^\infty \int_0^1 t_2 n (ax^n)^{k+1} dt_2 \\
&= -(2\pi i)^2 \sum_{k=0}^\infty n a^k \varepsilon^{nk} \frac{1}{2\pi i n k} \\
&= 2\pi i \log(1-a\varepsilon^n)
\end{aligned}$$

**Lemma 3.5.** *The integral*

$$\int_{\sigma_\varepsilon} \frac{d(1-ax^j)}{1-ax^j} \circ \frac{d(1-bx^k)}{1-bx^k}$$

*is 0 is $j \cdot k > 0$; otherwise it is equal to*

$$2\pi i sgn(k) d \sum_{l=1}^\infty \frac{a^{kl/d} b^{jl/d}}{l} = 2\pi i sgn(j) d \log(1-a^{k/d}b^{j/d}),$$

*where $sgn(k)$ is 1 if k is positive, -1 if k is negative; d is the greatest common divisor of j and k.*

*Proof.* Take the parametrization of $\sigma_\varepsilon : x \to \varepsilon e^{2\pi i t}, 0 \leq t \leq 1$. Then

$$\begin{aligned}
\frac{d(1-ax^j)}{1-ax^j} &= -(\sum_{n_1=0}^{+\infty} (ax^j)^{n_1}) d(ax^j) \\
&= -(\sum_{n_1=0}^{+\infty} j(ax^j)^{n_1+1}) \frac{dx}{x} \\
&= -2\pi i (\sum_{n_1=1}^{+\infty} j(a\varepsilon^j)^{n_1}) e^{2\pi i j n_1 t} dt
\end{aligned}$$

Hence

$$\int_{\sigma_\varepsilon} \frac{d(1-ax^j)}{1-ax^j} \circ \frac{d(1-bx^k)}{1-bx^k}$$

$$= \int_0^1 \left( \int_0^{t_2} -2\pi i (\sum_{n_1=1}^{+\infty} j(a\varepsilon^j)^{n_1}) e^{2\pi i j n_1 t_1} dt_1 \right) (-2\pi i (\sum_{n_2=1}^{+\infty} k(b\varepsilon^j)^{n_2}) e^{2\pi i j n t_2}) dt_2$$

$$= (2\pi i)^2 j \sum_{n_1=1}^{+\infty} \sum_{n_2=1}^{+\infty} jk a^{n_1} b^{n_2} \varepsilon^{jn_1+kn_2} \int_0^1 \left( \int_0^{t_2} e^{2\pi i j n_1 t_1} dt_1 \right) e^{2\pi i k n_2 t_2} dt_2$$

$$= 2\pi i \sum_{n_1=1}^{+\infty} \sum_{n_2=1}^{+\infty} \frac{k a^{n_1} b^{n_2} \varepsilon^{jn_1+kn_2}}{n_1} \int_0^1 (e^{2\pi i j n_1 t_2} - 1) e^{2\pi i k n_2 t_2} dt_2$$

$$= 2\pi i \sum_{n_1=1}^{+\infty} \sum_{n_2=1}^{+\infty} \frac{k a^{n_1} b^{n_2} \varepsilon^{jn_1+kn_2}}{n_1} \int_0^1 e^{2\pi i (jn_1+kn_2) t_2} dt_2$$

So the integral is 0 if $jn_1 + kn_2 \neq 0$. If $jn_1 + kn_2 = 0$, since $n_1$ and $n_2$ are positive integers, we have $n_1 = n|k|/\gcd(j,k)$ and $n_2 = n|j|/\gcd(j,k)$, where $n$ can be any positive integer. Write $d = \gcd(j,k)$, we have

$$\int_{\sigma_\varepsilon} \frac{d(1-ax^j)}{1-ax^j} \circ \frac{d(1-bx^k)}{1-bx^k}$$

$$= 2\pi i \sum_{n=1}^{+\infty} \frac{k a^{n_1} b^{n_2}}{n_1} = 2\pi i \sum_{n=1}^{+\infty} \frac{kd(a^{|k|/d})^n (b^{|j|/d})^n}{|k|n}$$

$$= -2\pi i \operatorname{sgn}(k) d \log(1 - a^{k/d} b^{j/d}) = 2\pi i \operatorname{sgn}(j) d \log(1 - a^{k/d} b^{j/d})$$

Combining all the above lemmas, we can easily see that

$$\frac{1}{2\pi i} \int_{\sigma_\varepsilon} \frac{df}{f} \circ \frac{dg}{g} = \sum_{j,k=1}^{\infty} (j,k) \log(1 - a_j^{k/(j,k)} b_{-k}^{j/(j,k)}) - \sum_{j,k=1}^{\infty} (j,k) \log(1 - a_{-j}^{k/(j,k)} b_k^{j/(j,k)})$$
$$+ \pi i \nu(f) \nu(g) + \nu(f) \log g(\varepsilon) - \nu(g) \log f(\varepsilon)$$

Now look at the loop $\sigma = \gamma_\varepsilon \sigma_\varepsilon \gamma_\varepsilon^{-1}$. By the composition of path formula 1.5, we have

$$\int_\sigma \frac{df}{f} \circ \frac{dg}{g} = \int_{\sigma_\varepsilon} \frac{df}{f} \circ \frac{dg}{g} + \int_{\gamma_\varepsilon} \frac{df}{f} \int_{\sigma_\varepsilon} \frac{dg}{g} + \int_{\sigma_\varepsilon} \frac{df}{f} \int_{\gamma_\varepsilon^{-1}} \frac{dg}{g}$$

**Lemma 3.6.** *Up to a multiple of $2\pi i$, we have*

$$\int_{\gamma_\varepsilon} \frac{df}{f} = \log f(\varepsilon) - \log f(P)$$

*Proof.* Write $\frac{df}{f}$ as

$$\nu(f)\frac{dx}{x} + \sum_{j=-N_f}^{-1} d\log(1 - a_j x^j) + d\log f^+.$$

Then it is easy to see that the lemma follows.

Now we can write down the result of the iterated integral of $\frac{df}{f}$ and $\frac{dg}{g}$ over the loop $\sigma$ using the composition of path 1.5 formula for iterated integrals and previous results

$$\begin{aligned}
&\frac{1}{2\pi i} \int_\sigma \frac{df}{f} \circ \frac{dg}{g} \\
&= \frac{1}{2\pi i} \left[ \left( \int_{\gamma_\varepsilon} \frac{df}{f} \circ \frac{dg}{g} + \int_{\gamma_\varepsilon^{-1}} \frac{df}{f} \circ \frac{dg}{g} + \int_{\gamma_\varepsilon} \frac{df}{f} \int_{\gamma_\varepsilon^{-1}} \frac{dg}{g} \right) \right. \\
&\quad \left. + \left( \int_\sigma \frac{df}{f} \circ \frac{dg}{g} + \int_{\gamma_\varepsilon} \frac{df}{f} \int_\sigma \frac{dg}{g} + \int_\sigma \frac{df}{f} \int_{\gamma_\varepsilon^{-1}} \frac{dg}{g} \right) \right] \\
&= \pi i \nu(f)\nu(g) + \sum_{j,k=1}^{\infty} (j,k)\log(1 - a_{-j}^{k/(j,k)} b_k^{j/(j,k)}) \\
&\quad - \sum_{j,k=1}^{\infty} (j,k)\log(1 - a_j^{k/(j,k)} b_{-k}^{j/(j,k)}) - \log \frac{f(P)^{\nu(g)} b_0^{\nu(f)}}{g(P)^{\nu(f)} a_0^{\nu(g)}}
\end{aligned}$$

The first part of the second line is equal to 0 by the product formula of iterated integrals.

The main theorem stated at the beginning of this section follows directly from this equation.

## 4. Applications

Using the reciprocity law provided in section 2, and the calculations given in section 4, it is easy to reprove the reciprocity law of the Contou-Carrère symbol:

**Theorem 4.1 (Reciprocity law of the Contou-Carrère symbol).** *Let $S$ be the set of zeroes and poles of $f$ and $g$. Then*

$$\prod_{s \in S} <f, g>_s = 1$$

*Proof:* Notice that $\int_{\sigma_i} \frac{df}{f} \int_{\sigma_j} \frac{dg}{g}$ is an integer multiple of $(2\pi i)^2$. For the part $\sum_{i=1}^{g} (\int_{\alpha_i} \frac{df}{f} \int_{\beta_i} \frac{dg}{g} - \int_{\beta_i} \frac{df}{f} \int_{\alpha_i} \frac{dg}{g})$, we can first write $f$ as $f_0 + \sum_{i=1}^{n} a_i f_i =$

$f_0(1+F)$, where $F = f_0^{-1} \sum_{i=1}^n a_i f_i$, then

$$\begin{aligned}
\int_{\alpha_i} \frac{df}{f} &= \int_{\alpha_i} \frac{df_0}{f_0} + \int_{\alpha_i} \frac{d(1+F)}{1+F} \\
&= \int_{\alpha_i} f^* \frac{dz}{z} + \int_{\alpha_i} \sum_{k=0}^N F^k d(1+F) \\
&= \int_{f_*\alpha_i} \frac{dz}{z} + \int_{\alpha_i} d \sum_{k=1}^N \frac{F^k}{k} \\
&= 2\pi i m
\end{aligned}$$

where $m$ is an integer. The sum on the second line is finite since $F$ is nilpotent. Hence $\sum_{i=1}^g (\int_{\alpha_i} \frac{df}{f} \int_{\beta_i} \frac{dg}{g} - \int_{\beta_i} \frac{df}{f} \int_{\alpha_i} \frac{dg}{g})$ is also an integer multiple of $(2\pi i)^2$. Combining all the results, divide both sides of the equation obtained in section 2, and then exponentiate, we get the reciprocity law of the Contou-Carrère symbol.

**Remark 4.2.** *If $f$ and $g$ are normal meromorphic functions, all the $a_j, b_k$ are 0 for $j < 0, k < 0$. Then we get back the Weil reciprocity.*

Using this new definition of the Contou-Carrère symbol, it is very easy to prove the following properties of the symbol:

**Proposition 4.3.**
  (1) $<f, gh> = <f,g><f,h>$
  (2) $<f, cf> = 1$
  (3) $<f, g^{-1}> = <f,g>^{-1} = <g, f>$
  (4) $<f, 1-f> = 1$, *if $f$ and $1-f$ are invertible*

*Proof:* (1) and (3) are easy to see by noticing the linearity of iterated integral and the facts that $\frac{d(gh)}{gh} = \frac{dg}{g} + \frac{dh}{h}$ and $\frac{dg^{-1}}{g^{-1}} = -\frac{dg}{g}$. (2) is obvious. To prove (4), first notice that

$$\int_\sigma \frac{df}{f} \circ \frac{d(1-f)}{1-f} = \int_\sigma f^\star(\frac{dz}{z}) \circ f^\star(\frac{d(1-z)}{1-z}) = \int_{f_\star\sigma} \frac{dz}{z} \circ \frac{d(1-z)}{1-z}$$

By the homotopy invariance of iterated integrals, this is equal to the integral over $\sigma_0^{n_1} \sigma_1^{n_2} \sigma_0^{n_3} \sigma_1^{n_4} \cdots \sigma_0^{n_{(N-1)}} \sigma_1^{n_N}$, where $\sigma_0$ is a simple loop around 0 starting and ending at $P$, and $\sigma_1$ is a simple loop around 1 starting and ending at $P$. Using the composition of path formula, it is easy to see that the result is an integer multiple of $(2\pi i)^2$, hence the resulting symbol is trivial.